\documentclass{article}
\usepackage[utf8]{inputenc}
\usepackage{amsthm}
\usepackage{amsmath}
\usepackage{amsfonts}
\usepackage{babel}
\usepackage{geometry}
\usepackage{fullpage}
\usepackage[colorlinks=false,hidelinks]{hyperref}
\usepackage{xurl}
\usepackage{color}
\usepackage{graphicx}
\usepackage[dvipsnames]{xcolor}
\usepackage[shortlabels]{enumitem}
\usepackage[nosort,capitalise]{cleveref}
\crefname{enumi}{}{}
\crefname{equation}{}{}
\usepackage{tikz}
\usetikzlibrary{matrix}
\usepackage{booktabs}

\definecolor{AfonsoBlue}{RGB}{30,65,123}

\makeatletter
\let\latex@fnsymbol\@fnsymbol
\renewcommand\@fnsymbol[1]{\ifcase#1\or*\or\S\or$\sharp$\or \P\else\@ctrerr\fi}
\newcommand{\restorefnsymbol}{\let\@fnsymbol\latex@fnsymbol}
\makeatother

\title{Randomstrasse101: Open Problems of 2024}

\author{
  Afonso S. Bandeira\thanks{ASB: Department of Mathematics, ETH Zürich, Rämistrasse 101, 8092 Zurich, Switzerland. \texttt{bandeira@math.ethz.ch}}
  \and
  Anastasia Kireeva\thanks{AK: Department of Mathematics, ETH Zürich, Rämistrasse 101, 8092 Zurich, Switzerland. \texttt{anastasia.kireeva@math.ethz.ch}}
  \and
  Antoine Maillard\thanks{AM: ARGO Team, Inria Paris, 48 Rue Barrault, 75013 Paris, France. \texttt{antoine.maillard@inria.fr}. Blog entry written while at the Department of Mathematics, ETH Zürich, Rämistrasse 101, 8092 Zurich, Switzerland.}
  \and
  Almut Rödder\thanks{AR: Department of Mathematics, ETH Zürich, Rämistrasse 101, 8092 Zurich, Switzerland. \texttt{almutmagdalena.roedder@ifor.math.ethz.ch}}
}

\date{\today}

\usepackage{natbib}
\usepackage{graphicx}

\newcommand{\RR}{\mathbb{R}}

\newcounter{sectionforenv}

\newtheorem{theorem}{Theorem}[sectionforenv]

\makeatletter
\@addtoreset{sectionforenv}{section}
\makeatother
\newtheorem{problem}{Problem}%[section]
\newtheorem{definition}[theorem]{Definition}

\newtheorem{conjecture}[problem]{Conjecture}
\newtheorem{openproblem}[problem]{Open Problem}

\newtheorem{remark}[theorem]{Remark}

\usepackage{tocloft}

\cftsetindents{section}{1em}{5em}
\begin{document}

\maketitle

\begin{abstract}
\texttt{Randomstrasse101} is a blog dedicated to Open Problems in Mathematics, with a focus on Probability Theory, Computation, Combinatorics, Statistics, and related topics. This manuscript serves as a stable record of the Open Problems posted in 2024, with the goal of easing academic referencing. The blog can currently be accessed at~\texttt{randomstrasse101.math.ethz.ch}.
\end{abstract}

%\addtocounter{section}{-1}

\section*{Introduction}

In this manuscript we include the blog entries in
\texttt{Randomstrasse101} of 2024, containing a total of fifteen Open Problems. \texttt{Randomstrasse101} is a blog created and maintained by our group at the Department of Mathematics at ETH Z\"{u}rich\footnote{The inspiration for the name of the blog will be clear to the reader after looking up the department's address and recalling the Probability Theory focus of the blog.}.

The focus is on mathematical open problems in Probability Theory, Computation, Combinatorics, Statistics, and related topics. The goal is not to necessarily write about the most important open questions in these fields but simply to discuss open questions and conjectures\footnote{Conjectures here should be interpreted as mathematical statements that we do not know not to be true, and for which a proof or a refutal would be interesting progress. We do not necessarily have very high confidence that conjectures here are true.} that each of us find particularly interesting or intriguing, somewhat in the same style as the first author's --now almost a decade old-- list of forty-two open problems~\cite{Afonso10L42P}. 

The blog was created and is currently maintained by Afonso S. Bandeira, Daniil Dmitriev, Anastasia Kireeva, Antoine Maillard, Chiara Meroni, Petar Nizic-Nikolac, Kevin Lucca, and Almut R\"{o}dder at ETH Z\"urich. Each blog entry is generally written by one author and the author list of this manuscript is the union of the entry authors in it.

Given the nature of this material, the writing of this manuscript is more informal than a typical academic text\footnote{If you would like to refer to an open question or blog entry, we encourage you to refer to this manuscript instead, since it is a more stable reference.}. Nevertheless, we hope it is useful and inspires thought on these questions. Happy solving!

\tableofcontents

\section{Welcome $\&$ Matrix Discrepancy (ASB)}

Let me start with one of my favourite Open Problems.

\begin{conjecture}[Matrix Spencer]\label{conj:1}
    There exists a positive universal constant $C$ such that, for all positive integers $n$, and all choices of $n$ self-adjoint $n\times n$ real matrices $A_1,\dots,A_n$ satisfying, for all $i\in[n]$, $\|A_i\|\leq 1$ (where $\|\cdot\|$ denotes the spectral norm) the following holds
    \begin{equation}\label{eq:matrixspencerbound}
    \min_{\varepsilon \in \{\pm1\}} \left\| \sum_{i=1}^n\varepsilon_iA_i\right\|\leq C\sqrt{n}.
    \end{equation}
\end{conjecture}

I first learned about this question from Nikhil Srivastava at ICM 2014, who pointed me to this very nice blog post of Meka~\cite{Meka-WindowsOnTheory-MSpencer}. To the best of my knowledge the question first appeared in a paper of Zouzias~\cite{Zouzias2011AMH}. I have also written about it in~\cite{Afonso10L42P} and more recently in~\cite{AfonsoOberwolfach10problems}.

The non-commutative Khintchine inequality (or a matrix concentration inequality) gives the bound up to a $\sqrt{\log n}$ factor.

In the particular case of commutative matrices, one can simultaneously diagonalize the matrices and the question reduces to a vector problem (representing the diagonal of the respective matrices) where the spectral norm is replaced by the $\ell_\infty$ norm. This is precisely the setting of Joel Spencer's celebrated ``Six Standard Deviations Suffice'' theorem~\cite{Spencer1985SixSD} which establishes the conjecture in that setting, for $C=6$. It is noteworthy that in the commutative setting a random choice does not give~\eqref{eq:matrixspencerbound}, the argument is of entropic nature. In a sense, in the other extreme situation in which the matrices behave ``very non-commutatively'' (or, more precisely, ``freely'') we know that a random choice of signs works, from our (myself, Boedihardjo, and van Handel) recent work on matrix concentration~\cite{BBvH-Free}. The fact that the reason for the conjecture to be true in these two extreme situations appears to be so different (one based on entropy the other on concentration), together with the fact that we still do not know if it is true in general makes this (in my opinion) a fascinating question!

In recent work, Hopkins, Raghavendra, and Shetty~\cite{Hopkins-MatrixSpencerQuantum} has established connections between this problem and Quantum Information, and Dadush, Jiang, and Reis~\cite{Dadush-MatrixSpencerMirror} establishes a connection with covering estimates for a certain notion of relative entropy.

More recently, Bansal, Jiang, and Meka~\cite{BJM-MatrixSpencer} has proved Conjecture~\ref{conj:1} for low-rank matrices (whose rank is $n$ over a polylog factor) using the matrix concentration inequalities in~\cite{BBvH-Free} within a very nice decomposition argument.

Motivated by Conjecture~\ref{conj:1} and the suspicion that ``ammount of commutativity'' may play an important role in this question, we (myself, Kunisky, Mixon, and Zeng)  proposed a Group theoretic special case of this conjecture in~\cite{BKMZ-Caley}.

\begin{conjecture}[Group Spencer]\label{conj:2}
  Let $G$ be a finite group of size $n$. Conjecture~\ref{conj:1} holds in the particular case in which $A_1,\dots,A_n$ are the $n\times n$ matrices corresponding to the regular representation of $G$.
\end{conjecture}

In~\cite{BKMZ-Caley} we prove Conjecture~\ref{conj:2} for simple groups, but the general case is still open. While the commutative case follows from Spencer's theorem, giving an explicit construction is not trivial (see this Mathoverflow post: \url{https://mathoverflow.net/q/441860}).

%Some of these problems are at the Oberwolfach Report \url{https://publications.mfo.de/handle/mfo/4149}

\section{Global Synchronization (ASB)}
%\addcontentsline{toc}{section}{Entry 2: Global Synchronization (ASB)}

In the 17th century, Christiaan Huygens (inventor of the pendulum clock) observed that pendulum clocks have a tendency to spontaneously synchronize when hung on the same board. This phenomenon of spontaneous synchronization of coupled oscillators has since become a central object of study in dynamical systems. 

We will focus here on spontaneous synchronization of $n$ coupled oscillators with pairwise connections given by a graph with adjacency matrix $A\in\RR^{n\times n}$. The celebrated Kuramoto model~\cite{Kuramoto1975} models the oscillators as gradient flow on $
\mathcal{E}(\theta) = \frac12\sum_{i,j=1}^n A_{ij}\left(1-\cos(\theta_i-\theta_j)\right),$
where the parameterization $\theta\in[0,2\pi[^n$ represents each oscillators angle in the moving frame of its natural frequency (for the sake of the sequel the derivation of this potential is not crucial you can see, e.g.~\cite{ABKSST-SynchExpander} for more details). This motivates the following definition which is central to what follows.

\begin{definition}
    We say an $n\times n$ matrix $A$ is globally synchronizing if the only local minima of $\mathcal{E}:\mathbb{S}^{n-1}\to \mathbb{R}$, parameterized by $\theta\in[0,2\pi[^n$ and given by
    \[
    \mathcal{E}(\theta) = \frac12\sum_{i,j=1}^n A_{ij}\left(1-\cos(\theta_i-\theta_j)\right),
    \]
    are the global minima corresponding to $\theta_i=c$, $\forall_i$. A graph $G$ is said to be globally synchronizing if its adjacency matrix is globally synchronizing.
\end{definition}

Abdalla, Bandeira, Kassabov, Souza, Strogatz, and Townsend~\cite{ABKSST-SynchExpander} recently showed that a random Erd\H{o}s-R\'enyi graph is globally synchronizing above the connectivity threshold with high probability (you can also see the Quanta article covering this at~\url{https://www.quantamagazine.org/new-proof-shows-that-expander-graphs-synchronize-20230724/}). The same paper also shows that a uniform $d$-regular graph is globally synchronizing with high probability for $d\geq 600$, but it leaves open the following question.

\begin{conjecture}[Globally Synchronizing Regular Graphs]\label{conj:8}
    % We say an $n\times n$ matrix $A$ is globally synchronizing if the only local minima of $\mathcal{E}:\mathbb{S}^{n-1}\to \mathbb{R}$, parameterized by $\theta\in[0,2\pi[^n$ and given by
    % \[
    % \mathcal{E}(\theta) = \frac12\sum_{i,j=1}^n A_{ij}\left(1-\cos(\theta_i-\theta_j)\right),
    % \]
    % are the global minima corresponding to $\theta_i=c$, $\forall_i$. A graph $G$ is said to be globally synchronizing if its adjacency matrix is globally synchronizing.
    A uniform random $3$-regular graph is globally synchronizing with high probability (probability going to $1$ as $n\to\infty$).
\end{conjecture}

The following question is motivated by trying to understand the Burer-Monteiro algorithmic approach to community detection in the stochastic block model~\cite{Bandeira-COLT16-BM}. The most recent progress is in~\cite{MABB-Z2Synch} which answers a high rank version of this question.

\begin{conjecture}[Global Synchrony with negative edges]\label{conj:9}
    Given any $\varepsilon>0$, the  $n\times n$ random matrix $A$ with zero diagonal and off-diagonal entries given by
    \[
    A_{ij} = 
    \left\{ \begin{array}{ccl}
    1 & \text{with probability} & \frac12+\delta\\
    -1 & \text{with probability} &
    \frac12-\delta,
    \end{array} \right.
    \]
    with $\delta \geq (1+\varepsilon)\sqrt{\frac{\log n}{2n}}$,
    is globally synchronizing with high probability.
\end{conjecture}

Another elusive question about global synchrony is the ``extremal combinatorics'' version of the question.

\begin{conjecture}[Density threshold for Global Synchrony]\label{conj:10}
For any $\varepsilon>0$, there exists $n>0$ and a graph $G$ on $n$ nodes such that the minimum degree of $G$ is at least $\left(\frac34-\varepsilon\right)n$ and $G$ is not globally synchronizing.
\end{conjecture}

Kassabov, Strogatz, and Townsend~\cite{KST-SynchDeg} showed an upper bound on the $\frac34$ threshold in Conjecture~\ref{conj:10}, and the same paper contains evidence that $\frac34$ is the correct threshold (by arguing that an upper bound below $\frac34$ would need to go beyond linear analysis of the dynamical system).

\section{Fitting ellipsoids to random points (AM)}
%\addcontentsline{toc}{section}{Entry 3: Fitting ellipsoids to random points (AM)}

In this entry I will discuss a seemingly simple question of high-dimensional stochastic geometry, which originated (as far as I know) in 
the series of works \cite{saunderson2011subspace,saunderson2012diagonal,saunderson2013diagonal}:
\begin{center}
\textit{For $n, d \gg 1$, when can a sequence of $n$ i.i.d.\ random points in $\mathbb{R}^d$, drawn from $\mathcal{N}(0, \mathrm{I}_d/d)$, be fitted by the boundary of a centered ellipsoid?}
\end{center}
\begin{figure}[!h]
    \centering
    \includegraphics[width=0.4\textwidth]{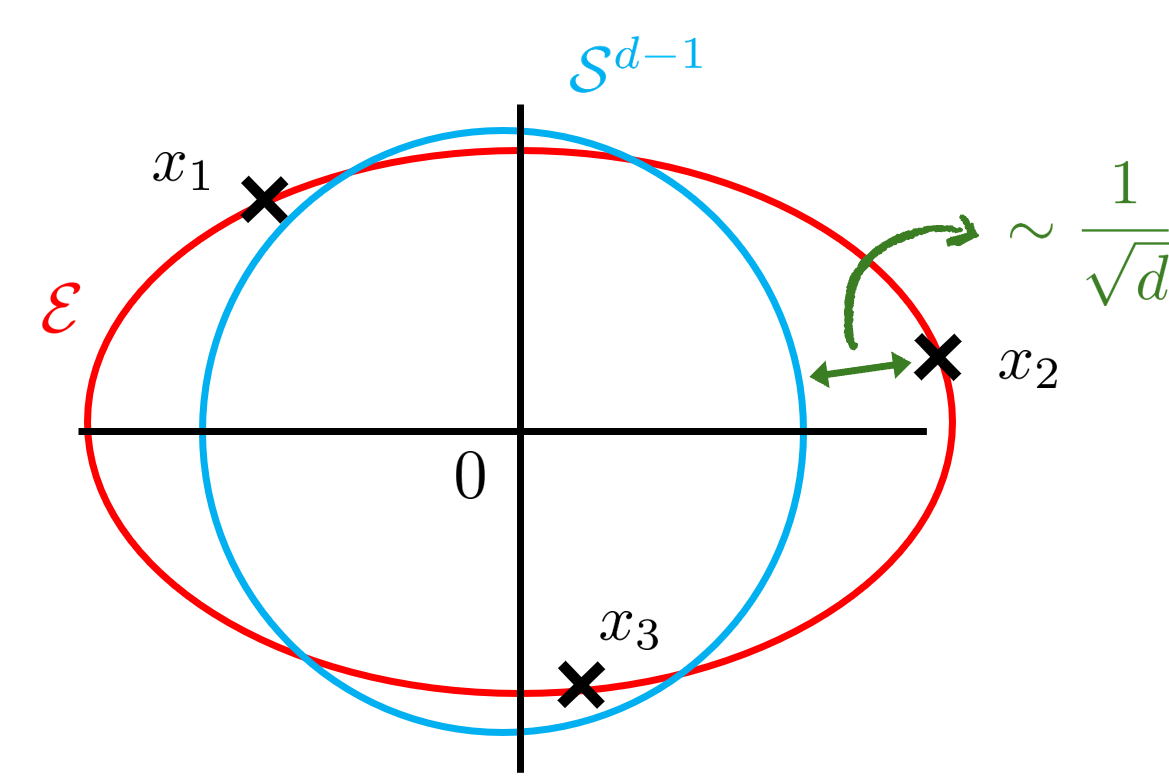}
    \caption{
        Fitting Gaussian random points $x_i \sim \mathcal{N}(0, \mathrm{I}_d/d)$ to an ellipsoid. Notice that the unit sphere itself is close to being a fit by simple concentration of measure: 
        a random $x \sim \mathcal{N}(0, \mathrm{I}_d/d)$ has (with high probability) distance $\mathcal{O}(1/\sqrt{d})$ to it.
    \label{fig:ellipsoid_fit}}
\end{figure}
The covariance being $\mathrm{I}_d / d$ is a convention which ensures $\mathbb{E}[\|x_i\|^2] = 1$: it is clear that assuming a generic positive-definite covariance $\Sigma \succ 0$ does not change this question, 
so we do not lose any generality with this assumption.

\medskip\noindent
The motivation of \cite{saunderson2011subspace,saunderson2012diagonal,saunderson2013diagonal} for this question came from statistics: they showed 
that the probability of existence of an ellipsoid fit was equal to the one of the success of a canonical convex relaxation (called Minimum-Trace Factor Analysis) in a problem of low-rank matrix decomposition.
Several other motivations were uncovered later on, and I recommend reading the introduction of \cite{potechin2023near} to learn more about them.

\medskip\noindent
Interestingly, it was soon conjectured that a sharp transition occurs in the regime $n/d^2 \to \alpha > 0$, exactly at $\alpha = 1/4$, see e.g.\ Conjecture~2 of \cite{saunderson2013diagonal}.
\begin{conjecture}[Ellipsoid fitting]\label{conj:ellipsoid_fitting}
    Let $n, d \geq 1$ and $x_1, \cdots, x_n \sim \mathcal{N}(0,\mathrm{I}_d/d)$. For any $\varepsilon > 0$: 
    \begin{align}
        \label{eq:ellipsoid_lower_bound}
        \limsup_{d \to \infty} \frac{n}{d^2} \leq \frac{1-\varepsilon}{4} \Rightarrow \lim_{d \to \infty} \mathbb{P}[\exists \mathcal{E} \textrm{ an ellipsoid fit to } (x_i)_{i=1}^n] = 1, \\
        \label{eq:ellipsoid_upper_bound}
        \liminf_{d \to \infty} \frac{n}{d^2} \geq \frac{1+\varepsilon}{4} \Rightarrow \lim_{d \to \infty} \mathbb{P}[\exists \mathcal{E} \textrm{ an ellipsoid fit to } (x_i)_{i=1}^n] = 0.
    \end{align}
\end{conjecture}

\medskip\noindent
\textbf{Upper bounds --}
Ellipsoid fitting can be written as a semidefinite program (SDP), since any origin-centered ellipsoid satisfies $\mathcal{E} = \{x \in \mathbb{R}^d \, : \, x^T S x = 1\}$ for some positive semidefinite symmetric matrix $S$, so:
\begin{equation}\label{eq:ellipsoid_sdp}
    \mathbb{P}[\exists \mathcal{E} \textrm{ an ellipsoid fit to } (x_i)_{i=1}^n] = 
    \mathbb{P}[\exists S \in \mathbb{R}^{d \times d} \, : \, S \succeq 0 \textrm{ and } x_i^T S x_i = 1 \textrm{ for all } i \in [n]].
\end{equation}
It's then not hard to convince oneself that $\{x_i^T S x_i = 1 \}_{i=1}^n$ is an independent system of $n$ linear equations on $S$: this already gives 
an upper bound of $\binom{d+1}{2} \simeq d^2/2$ to the number of points that can admit an ellipsoid fit (with high probability).
As of now, this ``silly'' argument (which does not take into account the constraint that $S \succeq 0$) is the best upper bound established for Conjecture~\ref{conj:ellipsoid_fitting}.

\medskip\noindent
\textbf{Lower bounds --}
On the other hand, there has been an abundance of works to establish lower bounds~\cite{saunderson2011subspace,saunderson2012diagonal,saunderson2013diagonal,kane2023nearly,potechin2023near,hsieh2023ellipsoid,tulsiani2023ellipsoid,bandeira2024fitting}.
Not diving into details, they essentially all rely on carefully picking a candidate solution $S^\star$ to the linear system $\{x_i^T S x_i = 1 \}_{i=1}^n$ (e.g.\ the least Frobenius norm solution), 
and establishing that $S^\star \succeq 0$ for small enough $n$. The best results in this vein are due to \cite{tulsiani2023ellipsoid,hsieh2023ellipsoid,bandeira2024fitting}, which independently proved the following.
\begin{theorem}[\cite{tulsiani2023ellipsoid,hsieh2023ellipsoid,bandeira2024fitting}]\label{thm:ellipsoid_lower_bound}
    Let $n, d \geq 1$ and $x_1, \cdots, x_n \sim \mathcal{N}(0,\mathrm{I}_d/d)$. There is $\delta > 0$ such that: 
    \begin{align*}
        \limsup_{d \to \infty} \frac{n}{d^2} \leq \delta \Rightarrow \lim_{d \to \infty} \mathbb{P}[\exists \mathcal{E} \textrm{ an ellipsoid fit to } (x_i)_{i=1}^n] = 1.
    \end{align*}
\end{theorem}
\noindent
While it is not known if these methods can be pushed all the way to $\delta = 1/4$,
it is conjectured that picking $S^\star$ as the minimal nuclear norm solution to the linear system $\{x_i^T S x_i = 1 \}_{i=1}^n$ gives this optimal value, and could be a path towards proving eq.~\eqref{eq:ellipsoid_lower_bound}~\cite{maillard2024fitting}.

\medskip\noindent
\textbf{The transition at $d^2/4$ --}
Notice that the $d^2/4$ threshold can be uncovered both by numerical simulations (solving a SDP can be done efficiently), but also by a heuristic argument: 
if the linear equations $\{\mathrm{Tr}[S x_i x_i^T] = 1 \}_{i=1}^n$ were replaced by $\{\mathrm{Tr}[S G_i] = 1 \}_{i=1}^n$ with $G_i$ Gaussian i.i.d.\ matrices, the existence of a 
sharp transition would be a direct consequence of Gordon's theorem~\cite{gordon1988milman}, and the value $d^2 / 4$ can even be recovered in this case as the squared Gaussian width of the cone of positive semidefinite matrices!
While this heuristic was well-known, we formalized it in \cite{maillard2023exact}, essentially by showing that the volume of the space of solutions is universal in both problems.
We established from it the following 
sharp transition result. 
\begin{theorem}[\cite{maillard2023exact}]
    \label{thm:ellipsoid_transition}
    For any $\varepsilon, M > 0$, we define 
    \begin{equation*}
        \mathrm{EFP}_{\varepsilon, M} \, : \, \exists S \in \mathbb{R}^{d \times d} \, : \, \mathrm{Sp}(S) \subseteq [0, M] \textrm{ and } \frac{1}{n} \sum_{i=1}^n |x_i^T S x_i - 1| \leq \frac{\varepsilon}{\sqrt{d}}.
    \end{equation*}
    Notice that the original ellipsoid fitting problem of Conjecture~\ref{conj:ellipsoid_fitting} is $\mathrm{EFP}_{0, \infty}$. We prove: 
    \begin{align}
       &\limsup_{d \to \infty} \frac{n}{d^2} = \alpha < \frac{1}{4} \, \Rightarrow  \,\exists M_\alpha > 0 \, : \, \forall \varepsilon > 0, \, \lim_{d \to \infty} \mathbb{P}[\mathrm{EFP}_{\varepsilon, M_\alpha}] = 1, \\
       &\liminf_{d \to \infty} \frac{n}{d^2} = \alpha > \frac{1}{4}  \,\Rightarrow  \,\exists \varepsilon_\alpha > 0 \, : \, \forall M > 0, \, \lim_{d \to \infty} \mathbb{P}[\mathrm{EFP}_{\varepsilon_\alpha, M}] = 0.
    \end{align}
\end{theorem}
\medskip\noindent
A few remarks are in order to clarify the conclusion of Theorem~\ref{thm:ellipsoid_transition}.

\medskip\noindent
\begin{itemize}
    \item The setting $\varepsilon \ll 1$ (but not going to $0$ with $d$) is precisely the regime where the problem is not trivially solved by the unit sphere (i.e.\ $S = \mathrm{I}_d$), since 
    $(1/n) \sum_{i=1}^n |x_i^T x_i - 1| = \mathcal{O}(1/\sqrt{d})$ (cf.\ also Fig.~\ref{fig:ellipsoid_fit}).
    \item In the regime $n / d^2 < 1/4$, Theorem~\ref{thm:ellipsoid_transition} shows that there exists ellipsoids which are: $(i)$ well-behaved (the spectral norm of $S$ is bounded), and $(ii)$ they fit the points $(x_i)$ up to an arbitrarily small error (but not going to $0$ with $d$).
    \item In the regime $n / d^2 > 1/4$, Theorem~\ref{thm:ellipsoid_transition} rules out any well-behaved ellipsoid (i.e.\ with bounded spectral norm) as a possible fit, 
    even allowing a small fitting error.
\end{itemize}

\medskip\noindent
Theorem~\ref{thm:ellipsoid_transition} provides the first mathematical result identifying a satisfiability transition in the ellipsoid fitting problem at the conjectured $d^2/4$ threshold. 
While we were not yet able to close Conjecture~\ref{conj:ellipsoid_fitting} from it, some parts seem tantalizingly close: e.g.\ the non-existence of ellipsoid fits for $n/d^2 > 1/4$ 
would now follow from showing that there is no ill-behaved ellipsoid fit (or that if there is an ill-behaved ellipsoid fit, there must also be a well-behaved one).

\medskip\noindent 
\textbf{Other works --}
To keep this post at a reasonable length, I have not mentioned some other recent works which are in some way related to this problem, and to which I contributed~\cite{maillard2024fitting,maillard2024bayes,erba2024bilinear}.
In \cite{maillard2024fitting} we use 
non-rigorous tools that originate in statistical physics to extend the conjecture to 
other classes of random points beyond Gaussian distributions (and the conjectured satisfiability threshold can then be very different from $d^2/4$!), but also to study the typical geometrical shape of ellipsoid fits, and even to predict analytically the performance of the methods used to prove the lower bounds discussed above (cf.\ the conjecture on the minimal nuclear norm solution mentioned below Theorem~\ref{thm:ellipsoid_lower_bound}).

\medskip\noindent
Finally, in \cite{maillard2024fitting} we realized that the ideas behind Theorem~\ref{thm:ellipsoid_transition} can be adapted to problems in statistical learning: precisely, we characterize optimal learning from data in a model of so-called ``extensive-width'' neural network, a regime which is both particularly relevant in practical applications and which had so far largely resisted to theoretical analysis. 
In another recent preprint~\cite{erba2024bilinear}, we build upon these ideas to study a toy model of a transformer architecture.

 \section{Detection in Multi-Frequency synchronization (AK)}
 %\addcontentsline{toc}{section}{Entry 4: Detection in Multi-Frequency synchronization (AK)}

In the study of average-case complexity, one is often interested in the critical threshold of the signal-to-noise ratio at which a problem becomes tractable. In this entry, we will consider a multi-frequency synchronization problem, where one obtains noisy measurements of the relative alignments of the signal elements through multiple frequency channels. We are interested whether it becomes possible to detect the signal using a time-efficient algorithm at a lower signal-to-noise ratio compared to the single-frequency model.

Let us define the synchronization problem more formally. In a general synchronization problem over a group $G$, one aims to recover a group-valued signal $g = (g_{1}, \dots, g_{n})\in G^n$  from its noisy pairwise information on $g_{k}g_{j}^{-1}$. One way to model such observations is through receiving a function of $g_{k}g_{j}^{-1}$ corrupted with additive Gaussian noise
$$
Y_{kj}= f(g_{k}g_{j}^{-1})+W_{kj}, \quad W_{kj}\sim\mathcal{N}(0, 1). 
$$
In this entry, we will focus on a setting where measurements are available for all pairs $(j,k)$.

Motivated by the Fourier decomposition of a non-linear objective of the non-unique games problem \cite{Bandeira_NonUniqueGames_2020}, we consider the model as receiving measurements through several frequency channels. For example, for $G = U(1) = \{e^{i\varphi}, \varphi \in [0, 2\pi)\}$, we consider measurements of the following form:
$$
\left\{\begin{aligned}
Y_1 & =\frac{\lambda}{n} x x^*+\frac{1}{\sqrt{n}} W^{(1)}, \\
Y_2 & =\frac{\lambda}{n} x^{(2)}\left(x^{(2)}\right)^*+\frac{1}{\sqrt{n}} W^{(2)}, \\
& \vdots \\
Y_L & =\frac{\lambda}{n} x^{(L)}\left(x^{(L)}\right)^*+\frac{1}{\sqrt{n}} W^{(L)},
\end{aligned}\right.
$$
where $W^{(1)},\dots,W^{(L)}$ are independent matrices, and $x^{(k)} = (x_1^k, \dots, x_n^k)$ denotes entrywise power. This case corresponds to the angular synchronization, where the objective is to determine phases $\varphi_1, \dots, \varphi_n \in [0, 2\pi]$ from their noisy relative observation $\varphi_k - \varphi_j \mod 2 \pi$, and equivalent to the synchronization over $SO(2)$, as $SO(2) \cong U(1)$.

In the general case of compact groups, we consider the Peter-Weyl decomposition instead. Informally, in this case, the Fourier modes correspond to irreducible representations, and each pairwise measurement corresponding to an irreducible representation $\rho$ is a block matrix with blocks given by
$$
Y_{k j}^\rho=\frac{\lambda}{n} \rho\left(g_k\right) \rho\left(g_j^{-1}\right)+\frac{1}{\sqrt{n}} W_{k j}^{(\rho)} \text { for each } \rho \in \Psi.
$$

The single-frequency model (i.e., when $L=1$) reduces to the Wigner spiked matrix model. In this model, the celebrated BBP transition \cite{BBP} postulates that above a critical value $\lambda \ge \lambda^* = 1$, detection is possible based on the top eigenvalue, while below this threshold, the top eigenvalue does not provide reliable information on the presence of the signal as $n$ grows to infinity. Moreover, for a variety of dense priors on signal $x$, this procedure is optimal in the sense that no algorithm can detect the signal reliably below the spectral threshold, $\lambda^*=1$, including algorithms with no constraints on the runtime \cite{perryOptimalitySuboptimalityPCA2016}. 

Does the situation change for the multi-frequency model? For simplicity, let us assume that the signal $x$ is sampled uniformly from $L$-th roots of unity, $x \sim \operatorname{Unif}(\{e^{2\pi i k / L}\}_{k=0}^{L-1})$. This corresponds to synchronization over the cyclic group $\mathbb Z / L \mathbb Z$. In this case, \cite{perryOptimalitySuboptimalityPCA2016} showed that the statistical threshold is $\lambda_{\text{stat}}=\Theta(\log L / L)$, i.e., below this threshold it is impossible to detect the signal, while above the threshold, there exists an inefficient algorithm for this task. The authors give upper and lower bounds on $\lambda_{\text{stat}}$ with exact constants, and the upper bound is lower than the spectral threshold for a sufficiently large number of frequencies, namely, $L\ge 11$. 

From the computational point of view, in \cite{kireeva2024computational}, it was shown that assuming the low-degree conjecture, no polynomial-time algorithm can detect the signal below the spectral threshold, $\lambda^*=1$, regardless of receiving additional information through additional channels. This result applies to the setting when the number of frequencies $L$ is constant compared to the dimension $n$ and when the signal is sampled uniformly from $SO(2)$ or any finite group of constant size. Combining this result with the optimality of PCA for the single-frequency model, this suggests that receiving only constant number of additional frequencies does not lower the computational threshold.
This opens up an intriguing question: how much more frequencies one requires so that the detection becomes possible using an efficient algorithm below the spectral threshold? We can formulate the following conjecture. 

\begin{openproblem}
    Consider synchronization model over $SO(2)$ or synchronization model over a finite group, where the signal $x$ is sampled uniformly over group elements. Find a scaling $L= L_n$ of number of frequencies such that there exists a polynomial-time algorithm that can detect the signal reliably for all $\lambda>\lambda_{\text{comp}, L}$ for some $\lambda_{\text{comp}, L}<1$ as $n\to\infty$.
\end{openproblem}

There is empirical evidence supporting the conjecture that, as $L$ diverges, the computational threshold becomes lower than the spectral threshold: numerical simulations in \cite{Gao2019MultiFrequencyPS} suggest that the variant of AMP with a carefully performed initialization can surpass the threshold when the dimension $n$ and the number of frequencies $L$ are comparable. 

%Not only is it interesting to determine the correct scaling of $L$ at which an efficient algorithm can detect the signal below the spectral threshold, but also to identify the correct computational threshold where the transition occurs.

In the computationally hard regime, we may also hope that there exists a sub-exponential algorithm whose runtime possibly depends on the signal strength. Such a tradeoff was observed in the sparse PCA problem \cite{ding2023subexponential}. For the synchronization model over finite groups,
\cite{perryOptimalitySuboptimalityPCA2016} proposed an algorithm for detecting the signal that works in the entire computationally hard regime. 
This algorithm involves maximizing a certain test function over all possible solutions $g \in G^n$ and thus has an exponential runtime. 

\begin{openproblem}
Consider synchronization model over $SO(2)$ or synchronization model over a finite group, where the signal $x$ is sampled uniformly over group elements. Fix number of frequencies $L$ that does not depend on $n$. In the computationally hard regime, i.e., when $\lambda \in (\lambda_{\text{stat}}, 1)$, does there exist a sub-exponential algorithm, i.e., an algorithm running in time $\exp{n^\delta}$ for some $\delta < 1$, that can detect the signal reliably as $n\to\infty$?
\end{openproblem}

%In \cite{kireeva2024computational}, we proved that in the considered setting, the problem is low-degree hard for degree $D = o(n^{1/3})$, which assuming the low-degree conjecture, suggests that the algorithm for detection should require at least $\exp(\tilde\Omega(n^{1/3}))$ runtime. It is interesting to see if this problem exhibits a smooth tradeoff between runtime and the signal strength analagous to the sparse PCA \cite{ding2023subexponential}

\section{Tensor PCA and the Kikuchi algorithm (ASB)}
%\addcontentsline{toc}{section}{Entry 5: Tensor PCA and the Kikuchi algorithm (ASB)}

In 2014, Andrea Montanari and Emile Richard~\cite{MontanariRichardTensorPCA} proposed a statistical model for understanding a variant of Principal Component Analysis in Tensor data. This is currently referred as the Tensor PCA problem. We will consider symmetric version of the problem in which the signal of interest is a point in the hypercube. Given $n,r$ and $\lambda$, the goal is to estimate (or detect) an unknown ``signal'' $x\in\{\pm1\}^n$ (drawn uniformly from the hypercube), from ``measurements'' as follows: for $i_1<i_2<...<i_r$, $$Y_{i_1,i_2,...,i_r} = \lambda x^{\otimes r} + Z_{i_1,i_2,...,i_r}$$ where $Z_{i_1,i_2,...,i_r}$ are iid $\mathcal{N}(0,1)$ random variables (and independent from $x$).

We will consider $r$ (and $\ell$, to be introduced below) fixed and $n\to\infty$, all big-O notation will be in terms of~$n$. Tensor PCA is believed to undergo a statistical-to-computaional gap: without regards for computational efficiency, it is possible to estimate $x$ for $\lambda=\Omega\left(n^{-\frac{r}{2}+\frac{1}{2}}\right)$. Efficient algorithms, such as the Sum-of-Squares hierarchy, are able to solve the problem at $\lambda=\tilde{\Omega}\left(n^{-\frac{r}{4}}\right)$, where $\tilde{\Omega}$ hides logarithmic factors. Local methods, such as gradient descent and approximate message passing succeed at $\lambda=\tilde{\Omega}\left(n^{-\frac12}\right)$. For $r=2$, the problem simply involves matrices, and indeed all these thresholds coincide. We point the reader to~\cite{WEM-Kikuchi} and references therein for more on each of these thresholds.

In 2019, Alex Wein, Ahmed El Alaoui, and Cris Moore~\cite{WEM-Kikuchi} proposed an algorithm for this problem based on the so-called Kikuchi Free Energy, it roughly corresponds to a hierarchy of message passing algorithms. They showed that this approach achieves (up to logarithmic factors) the threshold for the sum-of-squares approach, shedding light on the gap between the message passing and sum-of-squares frameworks.

We briefly describe the Kikuchi approach. We will focus on even $r$. There is a design parameter $\ell$ (with $n\gg \ell\geq \frac{r}2$) which we will consider fixed. The Kikuchi matrix $M$ is the $\binom{n}{\ell} \times \binom{n}{\ell}$ matrix (with rows and columns indexed by subsets $I\subset[n]$ of size $\ell$) given by
    \[
    M(\lambda)_{I,J} = \left\{ \begin{array}{ccl}
    Y_{I \Delta J} & \text{if} & |I \Delta J| = r,\\
    0 &\text{otherwise,}&
    \end{array} \right.
    \]
    where $I\Delta J = (I\cup J) \setminus (I\cap J)$ denotes the symmetric difference.
    
    The goal is to understand when the top of the spectrum of $M$ reveals the spike $x$. It is shown in~\cite{WEM-Kikuchi} that this happens for $\lambda = \tilde{\Omega}\left( n^{-
    \frac{r}4} \right)$. Since $Z$ is rotation invariant, we can assume WLOG that $x=\mathbf{1}$, the all-ones vectors. %Also we do the change of variables $\lambda = n^{\frac{r}4}\lambda$. 
    The following conjecture also appears in~
    \cite{AfonsoOberwolfach10problems} and is a reformulation of a conjecture in~\cite{WEM-Kikuchi}:

\begin{conjecture}[Kikuchi Spectral Threshold]
    Given $r,\ell, n$ positive integers satisfying $n\gg \ell\geq \frac{r}2$ and r even ($r$ and $\ell$ will be fixed and $n\to\infty$). For each $S\subset [n]$ with $|S|=r$ let $Z_S\sim\mathcal{N}(0,1)$, and all independent. Let $\lambda\geq 0$ and let $M$ be the $\binom{n}{\ell} \times \binom{n}{\ell}$ matrix (with rows and columns indexed by subsets $I\subset[n]$ of size $\ell$) given by
    \[
    M(\lambda)_{I,J} = \left\{ \begin{array}{ccl}
    \lambda + Z_{I \Delta J} & \text{if} & |I \Delta J| = r,\\
    0 &\text{otherwise,}&
    \end{array} \right.
    \]
    where $I\Delta J = (I\cup J) \setminus (I\cap J)$ denotes the symmetric difference.

    Let $\lambda^\natural_{r,\ell}$ denote the threshold at which eigenvalues ``pop-out'' of the spectrum of $M(\lambda)$: in other words $\lambda^\natural_{r,\ell}$ is the real number such that, for all 
    $\lambda>\lambda^\natural_{r,\ell}$, there exists $\varepsilon>0$ such that \(\mathbb{E}\lambda_{\max} M(\lambda)>(1+\varepsilon+o(1))\mathbb{E}\lambda_{\max} M(0),\)
    where $o(1)$ is a term that goes to zero as $n\to\infty$.

    For fixed $r$, we have $$n^{\frac{r}4}\lambda^\natural_{r,\ell}\to 0$$ as $\ell\to\infty$ (%more precisely $\limsup_{\ell\to\infty}\limsup_{n\to\infty}\lambda^\natural_{r,\ell}= 0$)
    note that this is after one has taken $n\to\infty$).
\end{conjecture}

This would establish the very interesting phenomenon that there is no sharp threshold for polynomial time algorithms in Tensor PCA in the following sense: This conjecture would imply that by increasing $\ell$ (while remaining $\Theta(1)$, corresponding to increasing the computational cost of the algorithm while keeping it polynomial time), one would be able to decrease the critical signal-to-noise ratio in the sense of $n^{\frac{r}4}\lambda^\natural_{r,\ell}$ getting arbitrarily close to zero. Since the bound in~\cite{WEM-Kikuchi} contains logarithmic factors on $n$ (present in $\tilde{\Omega}(\cdot))$, it does not allow to see this phenomenon for $\ell=\Theta(1)$.

For even $r$ and $\frac12r\leq\ell<\frac34r$ the threshold has been characterized in my work with Giorgio Cipolloni, Dominik Schr\"oder, and Ramon van Handel~\cite{BCSvH-Free2}, which is a follow up to the matrix concentration inequalities (leveraging intrinsic freeness) in my work with March Boedihardjo, and Ramon van Handel~\cite{BBvH-Free}.

\section{Did just a couple of deviations suffice all along? (ASB)}
%\addcontentsline{toc}{section}{Entry 6: Did just a couple of deviations suffice all along? (ASB)}

In September 2024, at a conference on Mathematical Aspects of Learning Theory in Barcelona (\url{https://www.crm.cat/mathematical-aspects-of-learning-theory/}) Dan Spielman gave a beautiful talk on discrepancy theory and some of its applications in clinical trials. Even though this was not the precise focus of the talk, it sparked a discussion that day about lower bounds for Joel Spencer's Six Deviations Suffice Theorem. I will describe some of the unanswered questions from this discussion. Many thanks to Dan Spielman, Tselil Schramm, Amir Yehudayoff, Petar Nizic-Nikolac, and Anastasia Kireeva with whom discussing this problem contributed to making the workshop a memorable week!

Given $n$, a positive integer, and $A$, an $n\times n$ matrix, we define the discrepancy of $A$, $\mathrm{disc}(A)$, as
\[
\mathrm{disc}(A) = \inf_{x\in\{\pm1\}^n}\|Ax\|_\infty.
\]
One of the motivations of this question is to understand how much smaller is $\inf_{x\in\{\pm1\}^n}\|Ax\|_\infty$ compared to the value, e.g. in expectation, when $x$ is drawn uniformly from the hypercube. We point to the references in this post for a more thorough discussion of the context, motivations, and state of the art.

In 1985, Joel Spencer~\cite{Spencer6Deviations} showed the cellebrated ``Six Deviations Suffice'' Theorem, which states that for any positive integer $n$ and any $A\in\{\pm1\}^{n\times n}$
 $\mathrm{disc}(A) \leq 6\sqrt{n}$. There have since been improvements on this constant, but the question we will pursue here concerns lower bounds.

When the condition of $\pm1$ entries is replaced by an $\ell_2$ condition on the columns, the corresponding question is the famous K\'{o}mlos Conjecture (which coincidentally is the first open problem in my lecture notes from around a decade ago~\cite{Afonso10L42P}).

\begin{conjecture}[K\'{o}mlos Conjecture]
There exists a universal constant $K$ such that, for all square matrices $A$ whose columns have unit $\ell_2$ norm, $\mathrm{disc}(A) \leq K$.
\end{conjecture}

It is easy to see that the statement of this conjecture implies Spencer's Theorem, since the columns of an $n\times n$, $\pm1$ matrix, have $\ell_2$ norm $\sqrt{n}$. The best current lower bound for K\'{o}mlos is a recent construction by Tim Kunisky~\cite{Tim-lowerboundKomlos} which shows a lower bound on $K$ of $1+\sqrt{2}$. Unfortunately, the resulting matrix is not a $\pm1$ matrix (even after scaling), and so it does not provide a lower bound for the Spencer setting.

\begin{openproblem}[How many deviations are needed afterall?]
What is the value of the following quantity?
 $$\sup_{n\in\mathbb{N}}\sup_{A\in\{\pm1\}^{n\times n}}\frac1{\sqrt{n}}\mathrm{disc}(A)$$
\end{openproblem}

It is easy to see that the matrix $\left[\begin{array}{cc}1 & 1 \\ 1 &-1\end{array}\right]$ has discrepancy $2=\sqrt{2}\sqrt{2}$, giving a lower bound of $\sqrt{2}$ to this quantity. Also, Dan Spielman mentioned a numerical construction (of a specific size) that achieved a slightly larger discrepancy value. To the best of my knowledge, there is no known lower bound above $2$, motivating the question: ``Did just a couple of deviations suffice all along?''.

Also, to the best of my knowledge, there is no known construction (or proof of existence) of an infinite family achieving a lower bound strictly larger than one. More explicitly, we formulate the following conjecture.
\begin{conjecture}\label{conjecture:limsupSpencerC} Given the definitions above, we have
 $$\limsup_{n\to\infty}\sup_{A\in\{\pm1\}^{n\times n}}\frac1{\sqrt{n}}\mathrm{disc}(A) > 1.$$
\end{conjecture}

\textbf{Hadamard matrices}

A natural place to search for a family proving a lower bound is to use the Sylvester construction of Hadarmard matrices. An Hadarmard matrix is a $\pm1$ matrix whose columns are all orthogonal (in other words, after scaling by $\frac{1}{\sqrt{n}}$ it is an orthogonal matrix). In 1867 James Joseph Sylvester proposed the following construction for sizes $n$ that are powers of $2$. $H_0 = 1$ and $H_k$ is the $2^{k}\times 2^{k}$ matrix given by $H_k = H_1 \otimes H_{k-1}$. In other words,
\[
H_k = \left[\begin{array}{cc}H_{k-1} & H_{k-1} \\ H_{k-1} &-H_{k-1}\end{array}\right].
\]
Recall that $n=2^k$. Since $\frac{1}{\sqrt{n}}H_k$ is an orthogonal matrix, it is easy to see that $\left\|H_kx\right\|_2=n$ for all $x\in\{\pm1\}^n$ and thus $\mathrm{disc}(H_k)\geq \sqrt{n}$. This is not tight for $H_1$ (as its discrepancy is $2$), however it is tight for $H_2$: e.g. taking the eigenvector $y = (1, 1, 1, -1)$.

To give an upper bound on $\mathrm{disc}(H_k)$ we proceed with an explicit construction of a $\pm1$ vector, using the (optimal) eigenvector $y$ for $H_2$. Let $x^\natural_{(0)}=1$, $x^\natural_{(1)}=(1, 1)$, and for $k \geq 2$, $x^\natural_{(k)}\in\{\pm1\}^{2^k}$ is given by
\[
x^\natural_{(k)} = y \otimes x^\natural_{(k-2)} =\left[\begin{array}{c}x^\natural_{(k-1)}  \\ x^\natural_{(k-1)}  \\x^\natural_{(k-1)}  \\ -x^\natural_{(k-1)}\end{array}\right].
\]
Since $H_k = H_1 \otimes (H_1 \otimes H_{k-2}) = H_2 \otimes H_{k-2}$, the mixed-product property gives
\[
H_kx^\natural_{(k)} = (H_2 \otimes H_{k-2}) (y \otimes x^\natural_{(k-2)}) = (H_2y) \otimes (H_{k-2}x^\natural_{(k-2)}) = (2 y) \otimes (H_{k-2}x^\natural_{(k-2)}).
\]
Inducting this argument shows that
\[
\left\|H_kx^\natural_{(k)}\right\|_\infty = \begin{cases}
\sqrt{2}\sqrt{2^k} \text{ if } k \text{ is odd,} \\
\sqrt{2^k} \text{ if } k \text{ is even.} 
\end{cases}
\]
While this settles the discrepancy of $H_k$ for $k$ even, the odd case is, to the best of my knowledge, open.

\begin{conjecture} Let $H_k$ be the $2^k\times 2^k$ Hadamard matrix obtained by the Sylvester construction (see above). For odd $k$, we have
\[
\mathrm{disc}(H_k) = \sqrt{2}\sqrt{2^k}.
\]
\end{conjecture}

Note that this would settle Conjecture~\ref{conjecture:limsupSpencerC}.

\begin{remark}[Boolean Fourier Analysis]
I am not going to discuss this at length here, but since $H_k$ essentially corresponds to the boolean Fourier transform for boolean functions on $k$ variables, the discrepancy of $H_k$ has a natural interpretation in Boolean Analysis. We point the reader to~\cite{BooleanAnalysisODonnell} for more on this subject.
\end{remark}

In a final note, I could not mention Hadamard matrices without stating a famous conjecture regarding their existence.

\begin{conjecture}[Hadamard Conjecture]
For any positive $n$ a multiple of $4$, there exists an $n\times n$ Hadamard matrix.
\end{conjecture}

\section{Sampling from the Sherrington-Kirkpatrick Model (AR)}
%\addcontentsline{toc}{section}{Entry 7: Sampling from the Sherrington-Kirkpatrick Model (AR)}

In recent years, significant progress has been made in understanding the mixing behavior of Glauber dynamics for Ising models, particularly in the Sherrington-Kirkpatrick (SK) model. Despite the classical Dobrushin-uniqueness condition being too restrictive for this model, recent breakthroughs employing techniques like stochastic localization and spectral analysis have improved the understanding of the fast-mixing regime. In this entry, we will discuss these improved mixing bounds for Glauber dynamics in the SK model and the open questions that remain.

Let us start with some preliminaries.
Consider the following measure on the hypercube $\{\pm 1\}^n$: 
$$
\mu_{J,h}(x) \propto \exp\left(\frac{1}{2} x^T J x +h^Tx\right),
$$
where $\propto$ corresponds to equal up to a normalizing constant.
Here, the interaction matrix $J \in \mathbb{R}^{n \times n}$ encodes pairwise interactions between spins, while $h$ represents the influence of an external field. We will particularly focus on $J = \beta W$ being a random matrix with $W \sim GOE(n)$ and inverse temperature $\beta>0$. This model was introduced by Sherrington and Kirkpatrick~\cite{sherrington1975solvable} and we refer to it as SK-model. The reader can focus on the case $h=0$. All statements below are with high probability on the disorder $W$.

Glauber Dynamics is an MCMC algorithm which chooses a site $i \in [n]$ uniformly at random and updates the spin according to $\mu$ conditioned on the remaining spins. The dynamics can be represented by the following transition kernel:
$$
P(x,y) = \begin{cases} \frac{1}{n}\mu(X_i = y_i \mid X_{-i} = x_{-i}) ,& \text{if } d_H(x,y) = 1\\
1-\frac{1}{n} \sum_{j=1}^n \mu(X_j \neq x_j \mid X_{-i} = x_{-i}) ,& \text{if } d_H(x,y) = 0\\
0, &\text{otherwise,}
\end{cases}
$$
representing the probability of transitioning from $x$ to $y$.
Here, $d_H$ denotes the Hamming distance in the hypercube.
Mixing in time $T_\epsilon$ means that, regardless of the starting point, the measure of the Markov chain at time $T_\epsilon$ is $\epsilon$-close to its stationary measure $\mu$ with respect to the total variation distance (although other notions of distances can also be used). The reader should think of $\epsilon$ as a small constant.

One celebrated and generally sufficient condition of Glauber Dynamics mixing in $O(n\log(n))$ time is the Dobrushin-uniqueness condition~\cite{dobruschin1968description} which requires $\max_{i}\sum_{j} |J_{ij}|<1$. While this condition is tight for some models (e.g. the Curie-Weiss-model), it only proves fast mixing for very high temperatures ($\beta = O(\frac{1}{\sqrt{n}})$) in the SK-model. In fact, a long standing open problem asked to extend the results above for some $\beta = O(1)$.

Indeed, in recent years, new proof techniques based on stochastic localization schemes have emerged and, when applied to the SK-model, significantly improved the bounds on $\beta$ for guaranteeing fast mixing. Eldan, Koehler and Zeitouni~\cite{eldan2022spectral} were able to prove the following general (deterministic) spectral condition: As long as
$$
\lambda_{\max}(J) -\lambda_{\min}(J) < 1-\delta
$$
for some $\delta >0$, then Glauber Dynamics mixes in $O(n^2)$ time with respect to total variation distance (big-O notation may hide constants depending on $\epsilon$). Applied to the SK-model this guarantees fast mixing for inverse temperature $\beta < \frac{1}{4}$.

The bound on the width of the spectrum of the Interaction Matrix proves to be tight again in the Curie-Weiss-Model. Moreover, based on the Low-degree-conjecture, Kunisky~\cite{kunisky2024optimality} gave evidence that no polynomial-time sampling algorithm exists that can succeed in general for all $J$ with $\lambda_{\max}(J) -\lambda_{\min}(J) \leq 1+\delta$. Hence, for improving bounds in the SK-model using spectral analysis, it does not suffice to measure the distance of the largest and the smallest eigenvalue of the interaction matrix, but additional properties of the interaction matrix in the SK model have to be used.
One nice example is the symmetry of the spectrum: In contrast to the all-ones-matrix which defines the Curie-Weiss-Interaction matrix, the spectrum of $W$ is almost symmetric around $0$. Anari, Koehler and Vuong~\cite{anari2024trickle} used exactly this property to improve the bound on $\beta$:
When
$$ \frac{|\lambda_{\max}(J)|}{|\lambda_{\min}(J)|} = 1 + o(1),
$$
Glauber Dynamics mixes with respect to Kullback-Leibler divergence in $O(n \log(n))$ time, when
$$
\lambda_{\max}(J) - \lambda_{\min}(J) < 1.18
$$
and $\lambda_{\min}(J) <0$.
The proof is based on a trickle down equation in a stochastic localization scheme and the authors introduce novel techniques to bound the operator norm of covariance matrices, depending on both the ratio of $\lambda_{\max}(J)$ and $\lambda_{\min}(J)$ and $\lambda_{\max}(J) - \lambda_{\min}(J)$.
This leads to the currently best-known upper bound $\beta < 0.295$ and it remains an open question whether this bound could be improved up to $1$:
\begin{openproblem}
     Consider the SK-model with an arbitrary external field. Does Glauber Dynamics fast mix (meaning in polynomial time) up to the replica symmetry breaking threshold $\beta_\ast =1$?
\end{openproblem}

El Alaoui, Montanari and Sellke~\cite{el2022sampling} proposed an algorithm for sampling from the SK-model without an external field that applies a stochastic localization scheme and makes use of the Approximate Message Passing framework. They proved fast mixing in Wasserstein Distance (which is a weaker metric than Total Variation Distance) up to $\beta <\frac{1}{2}$. Celentano~\cite{celentano2024sudakov} extended the argument to all $\beta$ up to the threshold $\beta_\ast =1$. 
Moreover, El Alaoui, Montanari and Sellke, proved in the same paper~\cite{el2022sampling} that it is impossible for stable algorithms to sample from the SK model without external field for  $\beta >1$, which also rules out fast mixing for their algorithm.

\section*{Updates}
\addcontentsline{toc}{section}{Updates}
While we generally do not expect to continuously update with progress on the questions discussed here, we try to include such updates if they happen (or we learn of them) before the posting of this manuscript. In this context, we are happy to report on two lines of progress: Conjecture 4 was answered by McRae in~\cite{mcrae2025benign}, following up on work of Rakoto Endor and Waldspurger~\cite{rakotoendor2024benign}. Furthermore, Rares-Darius Buhai~\cite{RaresDariusPrivate} showed that Conjecture 13 is false by identifying a connection between discrepancy of $H_k$ and a theory of non-linearity in Boolean Functions dating back, at least, to the work of Paterson and Wiedemann~\cite{PattersonWiedemann1983}; unfortunately, to the best of our knowledge, it does not rule out that the $H_k$ family of matrices can be used to establish Conjecture 12.

\bibliographystyle{alpha}
\bibliography{Arxiv2024}

\end{document}